\documentclass[11pt]{amsart}
\usepackage{mathtools}
\usepackage{graphicx}
\usepackage{color}
\usepackage[update,prepend]{epstopdf}
\usepackage{amsmath,amsfonts}
\usepackage{times}
\usepackage{cancel}
\usepackage[ansinew]{inputenc}
\usepackage{amssymb}
\usepackage{latexsym}
\usepackage{ulem}
\usepackage{cancel}
\newtheorem{rem}{Remark}[section]
\newcommand{\re}[1]{(\ref{#1})}

\numberwithin{equation}{section}

\def\openbox{$\sqcup\llap{$\sqcap$}$}
   
\let\qed=\endproof

\newtheorem{lem}{Lemma}[section]

\newtheorem{thm}{Theorem}[section]

\let\cal=\mathcal

\let \v=\varepsilon   
\let \a=\alpha   
 
\let\nt=\noindent
\let\ta=\theta
\let\bb=\mathbb

\let\re=\ref
\let\eqre=\eqref

  \def\openbox{$\sqcup\llap{$\sqcap$}$}
   
   \let\qed=\endproof

\begin{document}
\date{}

\title[Regularity ... semigroup ... interacting elastic systems]{Regularity and stability of the semigroup associated with some interacting elastic systems I: A degenerate damping case}
\author[Ka\"is Ammari]{Ka\"is Ammari}
\address{UR Analysis and Control of PDEs, UR 13ES64, Department of Mathematics, Faculty of Sciences of Monastir, University of Monastir, 5019 Monastir, Tunisia}
\email{kais.ammari@fsm.rnu.tn}

\author[Farhat Shel]{Farhat Shel}
\address{UR Analysis and Control of PDEs, UR 13ES64, Department of Mathematics, Faculty of Sciences of Monastir, University of Monastir, 5019 Monastir, Tunisia}
\email{farhat.shel@ipeit.rnu.tn}

\author[Louis Tebou]{Louis Tebou}
\address{Department of Mathematics and Statistics, Florida International University,
Modesto Maidique Campus, Miami, Florida 33199, USA}
\email{teboul@fiu.edu}
 
\begin{abstract}
In this paper, we examine regularity and stability issues for two damped abstract elastic 
 systems. The damping involves the average velocity and a fractional power $\theta$, with $\theta$ in $[-1,1]$, of the principal operator. The matrix operator defining the damping mechanism for the coupled system is degenerate. First, we prove that for $\theta$ in $(1/2,1]$, the underlying semigroup is not analytic, but is differentiable for $\theta$ in $(0,1)$; this is in sharp contrast with known results for a single similarly damped elastic system, where the semigroup is analytic for $\theta$ in $[1/2,1]$; this shows that the degeneracy dominates the dynamics of the interacting systems, preventing analyticity in that range.  Next, we show that for $\theta$ in $(0,1/2]$, the semigroup is of  certain Gevrey classes. Finally, we show that the semigroup decays exponentially for $\theta$ in $[0,1]$, and  polynomially for $\theta$ in $[-1,0)$. To prove our results, we use the frequency domain method, which relies on resolvent estimates. Optimality of our resolvent estimates is also established. Two examples of application are provided.
\end{abstract}


\maketitle
\tableofcontents

\section{Problem formulation and statements of main results} 
Let $H$ be a Hilbert space with inner product $(.,.)$ and norm $|.|$. Let $A$  be a positive unbounded self-adjoint operator, with domain $D(A)$ dense in the Hilbert space $H$. \\
 Set $V=D(A^\frac{1}{2})$. We assume that  $V\hookrightarrow H\hookrightarrow V'$, each injection being dense and compact, where $V'$ denotes the topological dual of $V$.\\
Let $a,~b$, and $\gamma$ be positive constants. Let $\theta\in[-1,1]$, and
consider the evolution system
\begin{equation}\label{e1}\begin{array}{lll}
&y_{tt}+aAy+\gamma A^{\theta}(y_t+z_t)=0\text{ in }(0,\infty)\\
&z_{tt}+bAz+\gamma A^{\theta}(y_t+z_t)=0\text{ in }(0,\infty)\\
&y(0)=y^0\in V,\quad y_t(0)=y^1\in H,\quad z(0)=z^0\in V,\quad z_t(0) =z^1 \in H.\end{array}\end{equation}
Introduce the Hilbert space ${\cal H}=V\times H\times V\times H$, over the field ${\bb C}$ of complex numbers, equipped with the norm
$$||Z||^2=a|A^\frac{1}{2}u|^2+|v|^2+b|A^\frac{1}{2}w|^2+|z|^2,\quad\forall Z=(u,v,w,z)\in{\cal H}.$$
Throughout this note, we shall assume:
\begin{equation}\label{xq1}
\exists a_0>0: |u|\leq a_0|A^\frac{1}{2}u|,\quad\forall u\in V.\end{equation}
Introduce the operator 
\begin{equation}\label{eac}
{\cal A}_\theta=\left( {\begin{array}{cccc}  
0&I&0&0\\  
-aA & -\gamma A^{\theta} &0&-\gamma A^{\theta}  \\
0&0&0&I\\
0 &-\gamma A^{\theta}&-bA & -\gamma A^{\theta}\\ 
\end{array}} \right)
\end{equation}
 with domain
$$D({\cal A}_\theta)=\Big\{(u,v,w,z)\in V^4;~ aAu+\gamma A^{\theta}(v+z)\in H,\text { and } bAw+\gamma A^{\theta}(v+z)\in H \Big\}.$$
One easily checks that for every $Z=(u,v,w,z)\in D({\cal A})$, 
\begin{equation}\label{e2}
\Re({\cal A}_\theta Z,Z)=- \gamma| A^{\frac{\theta}{2}}(v+z)|^2\leq0.\end{equation}so that the operator ${\cal A}_\theta$ is dissipative. Further, the operator ${\cal  A}_\theta$ is densely defined, so ${\cal A}_\theta$ is closable on ${\cal H}$. Therefore, the Lumer-Phillips Theorem shows that the operator ${\cal A}$ generates a strongly continuous semigroup of contractions $(S_\theta(t))_{t\geq0}$ on the Hilbert space ${\cal H}$. One also checks that 
$$i\bb R\subset\rho({\cal A}_\theta)$$where $\rho({\cal A}_\theta)$ denotes the resolvent set of ${\cal A}_\theta.$ This shows that the semigroup  $(S_\theta(t))_{t\geq0}$ is strongly stable on  the Hilbert space ${\cal H}$, thanks to the strong stability criterion of e.g. Arendt and Batty \cite{arb}. \\

This work was inspired by those of Chen and Triggiani \cite{ctra,ctrg}, where the authors considered the following abstract elastic system  
$$
\frac{d}{dt}\left( {\begin{array}{c}  
u\\  
v
\end{array}} \right)= \left( {\begin{array}{cc}  
0&I\\  
-A & -B_{\alpha}
\end{array}} \right)\left( {\begin{array}{c}  
u\\  
v
\end{array}} \right)
$$
on $D(A^{1\over2})\times H$. The operators $A$ and $B_\alpha$ are positive operators on the Hilbert space $H$ satisfying \eqre{xq1}, and are equivalent in a certain sense. They proved some regularity results for such systems. Similar results were proved by considering thermoelastic systems, see for example \cite{hly} where the authors present a complete  regularity and stability analysis  of the abstract system
$$
\begin{array}{lll}
&y_{tt}=-Ay+\delta A^{\alpha}z\text{ in }(0,\infty)\\
&z_{t}=-\delta A^{\alpha}y_t-\kappa A^{\beta}\text{ in }(0,\infty)\\
&y(0)=y^0,\quad y_t(0)=y^1,\quad z(0)=z^0.\end{array}
$$
The operator $A$ is as above, $\delta \neq 0$, $\kappa>0$ are real numbers, and $\alpha, \,\beta\in[0,1]$.

Our purpose in this work is twofold:
\begin{itemize}
\item  For $\theta\in(0,1]$, what type of regularity should we expect for the semigroup $(S_\theta(t))_{t\geq0}$? In other words, is the semigroup analytic for certain values of $\theta$ or not? If the semigroup is not analytic, is it of a certain Gevrey class? The reader should keep in mind that, in accordance with regularity results for a single elastic system damped as above, regularity is expected only for $\theta\in(0,1]$, e.g. \cite{ctra, ctrg, taylor}. In particular, responding to two conjectures of G. Chen and Russell \cite{cru}, S.P. Chen and Triggiani \cite{ctra,ctrg} proved for a single elastic system with a more general damping operator, which is equivalent to a fractional power $\alpha$ of the principal operator, that the semigroup is analytic for $\alpha$ in $[1/2,1]$, and of Gevrey class $s>1/2\alpha$ for $\alpha$ in $(0,1/2)$.
\item For $\theta\in[-1,1]$, what type of stability should we expect for the semigroup $(S_\theta(t))_{t\geq0}$?
\end{itemize}Thus, the main objective of our work is to analyze how the interplay of the dynamics of the coupled system affects the regularity or stability of the 
underlying semigroup.\vskip.2cm\noindent
Our findings are summarized in the following results:

\begin{thm}\label{neg} Assume that the constants $a$ and $b$ are distinct. \\
(i) For every $\theta\in(1/2,1]$, the semigroup $(S_\theta(t))_{t\geq0}$ is not analytic. In particular, for  every $\theta\in(1/2,1]$, and  every $r$ in $(2(1-\theta),1]$, we have :
\begin{align}\label{noan}
\limsup_{|\lambda|\to\infty} |\lambda|^{r}\left\Vert(i\lambda I - {\cal A}_\theta)^{-1}\right\Vert_{{\cal L}({\cal H})}=\infty.\end{align}
(ii) For every $\theta\in(0,1)$, the semigroup $(S_\theta(t))_{t\geq0}$ is differentiable, namely, there exists a positive constant $C$ such that we have the resolvent estimate:
\begin{align}\label{dif}
\log(|\lambda|)\left\Vert(i\lambda I-{\cal A}_{\ta})^{-1}\right\Vert_{{\cal L}({\cal H})}\leq C,\quad\forall\;\lambda\in{\mathbb R} \text{ with } |\lambda|\geq\lambda_0.
\end{align}for some large enough $\lambda_0>0$. The differentiability of the semigroup is valid for all $t>3K_0$, where \cite[p. 57]{pa}: $$K_0:=\limsup_{|\lambda|\to \infty}\left[\log(|\lambda|)\left\Vert(i\lambda I-{\cal A}_{\ta})^{-1}\right\Vert_{{\cal L}({\cal H})}\right].$$
\end{thm} 

\begin{thm}\label{reg}Assume that the constants $a$ and $b$ are distinct.
For every $\ta\in(0,1/2]$, the semigroup $(S_{\ta}(t))_{t\geq0}$ is of Gevrey class $\delta$ for every 
$\delta>1/s$, with $s=2\theta$ if $\theta\in(0,1/4]$, and $s=3\theta/(1+2\theta)$ if $\theta\in(1/4,1/2]$. More precisely, there exists a positive constant $C$ such that we have the resolvent estimate:
\begin{align}\label{geve}
|\lambda|^s\left\Vert(i\lambda I-{\cal A}_{\ta})^{-1}\right\Vert_{{\cal L}({\cal H})}\leq C,\quad\forall\;\lambda\in{\mathbb R}.
\end{align}Moreover, for every  $\ta\in(0,1/2)$ and every 
$r\in (2\theta,1]$, we have:
\begin{align}\label{optgev}
\limsup_{|\lambda|\to\infty} |\lambda|^{r}\left\Vert(i\lambda I - {\cal A}_\theta)^{-1}\right\Vert_{{\cal L}({\cal H})}=\infty.
\end{align}

\end{thm} 

\begin{thm}\label{stab} Assume that the constants $a$ and $b$ are distinct.
\begin{enumerate}
\item For every $\ta\in[0,1]$, the semigroup
$(S_{\ta}(t))_{t\geq0}$ is exponentially stable.  More precisely, there exist positive constants $K$ and $\xi$ such that:
\begin{align}\label{ese}
||S_{\ta}(t)||_{{\cal L}({\cal H})}\leq K\exp(-\xi t),\quad\forall\;t\ge 0.
\end{align}

\item For every $\ta\in[-1,0)$, the semigroup
$(S_{\ta}(t))_{t\geq0}$ is not exponentially stable. However, the semigroup is polynomially  stable, so that there exists a positive constant $K_0$ with:
\begin{align}\label{pse}
||S_{\ta}(t)Z^0||_{{\cal L}({\cal H})}\leq \frac{K_0 ||Z^0||_{D({\cal A}_{\ta})}}{(1+t)^{\frac{-1}{2\theta}}},\quad\forall Z^0\in D({\cal A}_{\ta}),\quad\forall\;t\ge 0.
\end{align}
Furthermore, the polynomial decay rate is optimal. More precisely, for  every  $\theta$ in $[-1,0)$ and every $r$ in $[0,~-2\theta)$, we have :
\begin{align*}
\limsup_{|\lambda|\to\infty} |\lambda|^{-r}\left\Vert(i\lambda I - {\cal A}_{\ta})^{-1}\right\Vert_{{\cal L}({\cal H})}=\infty.
\end{align*}
\end{enumerate}
\end{thm}

\begin{rem}One readily checks that if $a=b$ and we set $q=y-z$, where the pair $(y,z)$ is the solution of System \eqre{e1}, then $q$ satisfies the system

\begin{equation}\label{za1}\begin{array}{lll}
&q_{tt}+aAq=0\text{ in }(0,\infty)\\
&q(0)=y^0-z^0\in V,\quad q_t(0)=y^1-z^1\in H.\end{array}\end{equation} System \eqre{za1} is a conservative system; no regularity should then be expected of System \eqre{e1}, since it can be decoupled into $q=y-z$ and $p=y+z$, where the $p$-system is damped as \eqre{e1}, and as such, enjoys the regularity demonstrated in \cite{ctra, ctrg, taylor}, the $q$-system being conservative enjoys no regularity whatsoever. If we further assume 
$y^0\not=z^0$ or $y^1\not=z^1$, then the energy of system \eqre{za1} is nonzero for all times; so no stability is to be expected of System \eqre{e1}. Indeed the synchronization phenomenom occurs; this is in agreement with the findings in e.g. \cite{litr1}.
\end{rem} 

\begin{rem}Concerning the resolvent estimate leading to the differentiability of the semigroup associated with System \eqre{e1}, our proof shows that we can have the stronger resolvent estimate:
\begin{align}\label{difa}
\forall r\geq1,\quad\exists \lambda_0=\lambda_0(r)>0: (\log(|\lambda|))^r\left\Vert(i\lambda I-{\cal A}_{\ta})^{-1}\right\Vert_{{\cal L}({\cal H})}\leq C,\quad\forall\;\lambda\in{\mathbb R} \text{ with } |\lambda|\geq\lambda_0.
\end{align}\end{rem}

\begin{rem} It is quite surprising that analyticity fails for the coupled elastic system for $\theta$ in $(1/2,1]$. In that same range, the semigroup is differentiable though, except for $\theta=1$. We remind the reader that for a single  similarly damped elastic system, analyticity holds for $\theta$ in $[1/2,1]$, \cite{ctra, ctrg}. However, one should keep in mind that in the present situation, we are dealing with a degenerate system, as the matrix defining the damping is degenerate; it then appears that for $\theta$ in $(1/2,1]$, the degeneracy dominates the dynamics of the system, and precludes analyticity of the underlying semigroup. We think that the situation would be different in the non-degenerate case, which we plan to analyze in a subsequent work. For $\theta$ in $(0,1/4]$, we obtain the same regularity result as for a single similarly damped elastic system. For $\theta$ in $(1/4,1/2]$, we obtain a weaker regularity result, and we think that in this range too, the degeneracy is at work, but its action is milder than in the range $(1/2,1]$. Note that  estimate \eqref{optgev} shows that the resolvent estimate established for $\theta$ in $(0,1/4]$ is optimal. That same estimate \eqref{optgev}, which is actually valid for all $\theta$ in $(0,1/2)$,  also suggests that the resolvent estimate obtained for $\theta$ in $(1/4,1/2]$ is not optimal. 

\medskip

As for our polynomial decay estimate, we point out that for $\theta$ in $[-1,0)$, no exponential decay is to be expected, as the damping plays the role of a compact perturbation of an otherwise conservative system; in this case a result of Gibson \cite{gib} or Triggiani \cite{trig} ensures that exponential decay fails. Anyway the optimality of our polynomial decay estimate also precludes the exponential decay of the semigroup. We also point out that we obtain the same polynomial decay rate as for a single similarly damped elastic system, e.g. \cite{lzh}.
\end{rem}

\section{Some technical Lemmas}

\begin{lem}\label{expdec}  {\bf (\cite{hf, pr})} {\sl
Let ${\cal A}$ be the generator of a bounded $C_0$ semigroup
$(S(t))_{t\geq0}$ on a Hilbert space ${\cal H}$. Then
$(S(t))_{t\geq0}$ is exponentially stable if and only if:\par {\tt
i)} $i{\bb R}\subset\rho({\cal A})$, and\par  {\tt ii)}
$\sup\{||(ib-{\cal A})^{-1}||;~b\in {\bb R}\}<\infty$, where
$\rho({\cal A})$ denotes the resolvent of} ${\cal A}$.\end{lem}

\begin{lem}\label{pdec} {\bf (\cite{bort})}  Let ${\cal A}$ be the
generator of a bounded $C_0$ semigroup $(S(t))_{t\geq0}$ on a
Hilbert space ${\cal H}$ such that $i{\bb R}\subset\rho({\cal A})$,
where $\rho({\cal A})$ denotes the resolvent of ${\cal A}$. Then
$(S(t))_{t\geq0}$ is polynomially stable, {\tt viz.}, there are
positive constants $M$ and $\a$ that are independent of the initial
data such that
\begin{equation*}||S(t)Z^0||_{{\cal H}}\leq {M||Z^0||_{D({\cal A})}\over
(1+t)^{1\over\a}},\quad\forall t\geq0,\quad Z^0\in D({\cal
A}).\end{equation*} if and only if \begin{equation*}\exists C_0>0: ||(ib-{\cal
A})^{-1}||_{{\cal L}({\cal H})}\leq C_0|b|^{\a},\forall b\in{\bb R}\hbox{ with }
|b|\geq1.\end{equation*}\end{lem} 
General and weaker versions of Lemma \ref{pdec} may be found, respectively, in
\cite{seifert} and \cite{beps, lyu, bdu, lir}.

\begin{lem}\label{difc} (\cite[Chap. 2, p. 57]{pa}) 
 Let $T=(T (t))_{t\ge 0}$ be a strongly continuous semigroup on a Hilbert space $X$, with $||T(t)\leq Me^{\omega t}$. Let $A$ denote the infinitesimal generator of $A$. If for some $\mu\geq\omega$:
 
\begin{align*}
\limsup_{|\lambda|\to\infty}\left[ \log(|\lambda|)||(\mu+i\lambda I - A)^{-1}||_{{\cal L}(X)}\right]=C<\infty.
\end{align*} 
Then $T=(T (t))_{t\ge 0}$ is differentiable for $t > 3C$.
\end{lem}

\begin{lem}\label{anal} (\cite[Chap. 1, p. 5]{liz}) 
 Let $T=(T (t))_{t\ge 0}$ be a strongly continuous semigroup of contractions on a Hilbert space $X$, with infinitesimal generator $A$. Suppose that 
 
 $$i{\bb R}\subset\rho(A),$$where $\rho(A)$ denotes the resolvent of $A$.\\ The semigroup $T=(T (t))_{t\ge 0}$ is analytic if and only if 
 
\begin{align*}
\limsup_{|\lambda|\to\infty} |\lambda|||(i\lambda I - A)^{-1}||_{{\cal L}(X)}<\infty.
\end{align*} 
\end{lem}

\begin{lem}\label{gevc} (\cite{taylor}) 
 Let $T=(T (t))_{t\ge 0}$ be a strongly continuous and bounded semigroup on a Hilbert space $X$.
 Suppose that the infinitesimal generator $A$ of the semigroup $T$ satisfies the following estimate, for some $0<\alpha<1$:
\begin{align}\label{ege}
\limsup_{|\lambda|\to\infty} |\lambda|^\alpha||(i\lambda I - A)^{-1}||_{{\cal L}(X)}< \infty.
\end{align} 
Then $T=(T (t))_{t\ge 0}$ is of Gevrey class $\delta$ for $t > 0$, for every $\delta > {1\over\alpha}$.
\end{lem}

\section{ Proof of Theorem \ref{neg}}
This proof is divided into two parts. In the first part, we prove the non-analyticity of the semigroup, while the second part is devoted to establishing the differentiability of the semigroup.\\
{\bf Part 1: The semigroup is not analytic for $\pmb{\theta\in(1/2,1]}$.}
Here, we proceed as in \cite{ktw,twh}. Let $\theta\in(1/2,1]$.
We are going to show that there exist a sequence of  positive real numbers  $(\lambda_n)_{n\geq1}$, and for each $n$, an element $Z_n\in {\cal D}({\cal A})$
such that for every $r\in(2(1-\theta),1]$, one has:
\begin{align}\label{eazc}
\lim_{n\to\infty}\lambda_n=\infty,\quad ||Z_n||=1, \quad \lim_{n\to\infty}\lambda_n^{-r}||(i\lambda_n-{\cal A}_\theta)Z_n||=0.
\end{align}
Indeed, if we have sequences $\lambda_n$ and $Z_n$ satisfying \eqref{eazc}, then we set
\begin{align}\label{eaze}
V_n=\lambda_n^{-r}(i\lambda_n-{\cal A}_\theta)Z_n,\qquad U_n=\frac{V_n}{||V_n||}.
\end{align}
 Therefore,
$||U_n||=1$ and
\begin{align}\label{eazf}
\lim_{n\to\infty}\lambda_n^r||(i\lambda_n-{\cal A}_\theta)^{-1}U_n||=\lim_{n\to\infty}\frac{1}{||V_n||}=\infty,
\end{align}
which would establish the claimed result, thereby completing the proof of Theorem \ref{neg}.
 Thus, it remains to prove the existence of such sequences.\\
 We shall borrow some ideas from \cite{twh}. For each
$n\geq1$, we introduce the eigenfunction $e_n$, with
$|e_n|=1, \, A e_n= \omega_n e_n,$ and where $(\omega_n)$ is an increasing sequence of positive real numbers with $\displaystyle \lim_{n\to\infty}\omega_n=\infty.$ We seek $Z_n$ in
the form $Z_n=(a_ne_n,i\lambda_na_ne_n,c_ne_n,i\lambda_nc_ne_n)$, with $\lambda_n$ and the complex numbers
$a_n$ and $c_n$ chosen such that $Z_n$ fulfills the
desired conditions.\\ For this purpose, set $\alpha=(a+b)/2$, $\beta=(a-b)/2$, and for every $n\geq1$, set:
\begin{align}\label{eazd}
\lambda_n=\sqrt{\alpha \omega_n}.
\end{align}
With that choice, we readily check that:
\begin{equation}\label{eqpq}\begin{array}{lll}
(i\lambda_n-{\cal A}_\theta)Z_n&= \begin{pmatrix}
   0\\
\left[(a\omega_n-\lambda_n^2)a_n+i\lambda_n\gamma\omega_n^\theta( a_n+ c_n)\right]e_n\\
0\\ (-\lambda_n^2+b\omega_n)c_n+i\lambda_n\gamma\omega_n^\theta( a_n+c_n))
e_n
 \end{pmatrix}\\
&= \begin{pmatrix}
   0\\
\left[(\beta+i\alpha^{1\over2}\gamma\omega_n^{\theta-{1\over2}}) a_n+i\alpha^{1\over2}\gamma\omega_n^{\theta-{1\over2}} c_n\right]\omega_ne_n\\
0\\ \left[(-\beta+i\alpha^{1\over2}\gamma\omega_n^{\theta-{1\over2}}) c_n+i\alpha^{1\over2}\gamma\omega_n^{\theta-{1\over2}} a_n\right]\omega_ne_n
 \end{pmatrix},\text{ by }\eqre{eazd}.\\
\end{array}\end{equation}
If for each $n\geq1$, we set
\begin{align}\label{lak1}
a_n=-{i\alpha^{1\over2}\gamma\omega_n^{\theta-{1\over2}}c_n\over\beta+i\alpha^{1\over2}\gamma\omega_n^{\theta-{1\over2}}},\text{ and }  \omega_n^{1\over2}c _n=\alpha_n+i\beta_n,\end{align}
where $\alpha_n$ and $\beta_n$ are real numbers to be specified later on. It then follows from \eqref{eqpq}:

 \begin{equation}\label{lak2}\begin{array}{lll}
(i\lambda_n-{\cal A})Z_n&= \begin{pmatrix}
   0\\
0\\0\\
 \left[\left(-\beta+i\alpha^{1\over2}\gamma\omega_n^{\theta-{1\over2}}+\frac{\alpha\gamma^2\omega_n^{2\theta-1}}{\beta+i\alpha^{1\over2}\gamma\omega_n^{\theta-{1\over2}}}\right) c_n\right]\omega_ne_n
 \end{pmatrix}\\
 &= \begin{pmatrix}
   0\\
0\\
0\\ -{\beta^2\over \beta+i\alpha^{1\over2}\gamma\omega_n^{\theta-{1\over2}}}\omega_n c_ne_n
 \end{pmatrix}\end{array}\end{equation}

Therefore, we have (keeping in mind that $\theta>1/2$):
\begin{align}\label{lak01}
&\lim_{n\to\infty}\lambda_n^{-2r}||(i\lambda_n-{\cal A}) Z_n||^2\notag\\
&=\lim_{n\to\infty}{\beta^4\over \beta^2+\alpha\gamma^2\omega_n^{2\theta-1}}|c_n|^2\lambda_n^{-2r}\omega_n^2|e_n|_2^2=\lim_{n\to\infty}{\alpha^{-r}\beta^4\omega_n^{1-r}\over \beta^2+\alpha\gamma^2\omega_n^{2\theta-1}}(\alpha_n^2+\beta_n^2)\\
&=\lim_{n\to\infty}{\alpha^{-r-1}\gamma^{-2}\beta^4\omega_n^{2-r-2\theta}}(\alpha_n^2+\beta_n^2)=0,\text{ as } r>2(1-\theta), \text{ and }\notag
\end{align}provided that the sequences $(\alpha_n)$ and $(\beta_n)$ are built in such a way that they converge to some suitable real numbers, and $||Z_n||=1$, for $n$ large enough.\\ One checks that
\begin{align}\label{eqpp}
||Z_n||^2&=a\omega_n|a_n|^2+\lambda_n^2|a_n|^2+b\omega_n|c_n|^2+\lambda_n^2|c_n|^2
\notag\\
&=(a+\alpha)\omega_n|a_n|^2+(b+\alpha)\omega_n|c_n|^2)\notag\\
&=\left[(a+\alpha){\alpha\gamma^2\omega_n^{2\theta-1}\over\beta^2+\alpha\gamma^2\omega_n^{2\theta-1}}+(b+\alpha)\right]\omega_n|c_n|^2\\
&=\left[(a+\alpha)\left(1-{\beta^2\over\beta^2+\alpha\gamma^2\omega_n^{2\theta-1}}\right)+(b+\alpha)\right]\omega_n|c_n|^2\notag\\
&=\left[2(a+b)-{(a+\alpha)\beta^2\over\beta^2+\alpha\gamma^2\omega_n^{2\theta-1}}\right](\alpha_n^2+\beta_n^2),\text{ as } 2\alpha=a+b.\notag\end{align}
Now, we shall build sequences $(\alpha_n)$ and $(\beta_n)$. For each $n\geq 1$, set
$$\alpha_n=(\alpha_0+r_n),\quad \beta_n=(\beta_0+r_n),$$where the real numbers $\alpha_0$ and $\beta_0$ are chosen with
\begin{equation}\label{gev1}2(a+b)(\alpha_0^2+\beta_0^2)=1,\end{equation}and $r_n$ is to be conveniently chosen in the sequel, such that $r_n\to0$ as $n\nearrow\infty$.\\
Notice that with those choices, we now have, setting, $\zeta_n=(a+\alpha)\beta^2/4\alpha[(b+\alpha)\beta^2+4\alpha^2\gamma^2\omega_n^{2\theta-1}]$:
\begin{equation}\label{gev2}\begin{array}{lll}
||Z_n||^2&=\left[2(a+b)-{(a+\alpha)\beta^2\over\beta^2+\alpha\gamma^2\omega_n^{2\theta-1}}\right](\alpha_n^2+\beta_n^2)\\
&=\left[2(a+b)-{(a+\alpha)\beta^2\over\beta^2+\alpha\gamma^2\omega_n^{2\theta-1}}\right](\alpha_0^2+\beta_0^2+2(\alpha_0+\beta_0)r_n+2r_n^2)\\
&=1-{(a+\alpha)\beta^2\over\beta^2+\alpha\gamma^2\omega_n^{2\theta-1}}(\alpha_0^2+\beta_0^2) \\
&\hskip.3cm+\left[2(a+b)-{(a+\alpha)\beta^2\over\beta^2+\alpha\gamma^2\omega_n^{2\theta-1}}\right](2(\alpha_0+\beta_0)r_n+2r_n^2)\\
&=1-{(a+\alpha)\beta^2\over4\alpha(\beta^2+\alpha\gamma^2\omega_n^{2\theta-1})} +\left[4\alpha-{(a+\alpha)\beta^2\over\beta^2+\alpha\gamma^2\omega_n^{2\theta-1}}\right](2(\alpha_0+\beta_0)r_n+2r_n^2)\\
&=1-{(a+\alpha)\beta^2\over4\alpha(\beta^2+\alpha\gamma^2\omega_n^{2\theta-1})} +{(b+\alpha)\beta^2+4\alpha^2\gamma^2\omega_n^{2\theta-1}\over\beta^2+\alpha\gamma^2\omega_n^{2\theta-1}}(2(\alpha_0+\beta_0)r_n+2r_n^2),\\
\end{array}
\end{equation} as $2(a+b)(\alpha_0^2+\beta_0^2)=1$, and $3\alpha-a=b+\alpha$.\\Thus, $||Z_n||=1$ if and only if the quadratic equation 
\begin{equation}\label{gev3}
2r_n^2+2(\alpha_0+\beta_0)r_n-\zeta_n=0,\end{equation}
has at least one real root. \\ Now, the discriminant of the quadratic equation is positive, and its roots are given by:
$$r_n^{\pm}=\frac{-(\alpha_0+\beta_0)\pm\sqrt{(\alpha_0+\beta_0)^2+2\zeta_n}}{2}.$$ Notice that $\zeta_n\to0$ as $n\nearrow\infty$. Accordingly, as $r_n$ must converge to zero, for $\alpha_0+\beta_0>0$, we choose $r_n=r_n^{+}$, and for $\alpha_0+\beta_0<0$, we  choose $r_n=r_n^{-}$. When $\alpha_0+\beta_0=0$, any of the two roots will do. This  completes the proof of the claimed lack of analyticity.\\
\\
{\bf Part 2: Differentiability for $\pmb{\theta\in(0,1)}$.} From now on, we shall use the following notations
$$||u||_{\tau}=|A^{\tau\over2}u|\quad\forall\tau\in[-1,1],\quad\forall u\in V.$$
Let $r>0$ be an arbitrary real number. Let $\lambda_0=\lambda_0(r)>0$ be so large that for every $\lambda\in\mathbb R\text{ with }|\lambda|\geq\lambda_0$: 
\begin{equation}\label{log}
\frac{\log(|\lambda|)}{|\lambda|^r}\leq1.\end{equation}
Thanks to Lemma \ref{difc}, (note that we may choose $\mu=0$ in our case),  it is enough to prove that there exists a positive constant $C_0$ such that:
\begin{align}\label{dif0}
 \sup\left\{\log(|\lambda|)||(i\lambda I-{\cal A}_{\ta})^{-1}||_{{\cal L}({\cal H})};~ \lambda\in \mathbb R\text{ with }|\lambda|\geq\lambda_0\right\}\leq C_0.
\end{align}
The constant $C_0$ may vary from line to line and depends on the parameters of the system, but not on the frequency variable $\lambda$.\\
 To prove \eqref{dif0}, we will
show that there exists $C_0>0$ such that for every $U\in{\cal
H}$, one has:
\begin{align}\label{dif1}
\log(|\lambda|)||(i\lambda I-{\cal A}_{\ta})^{-1}U||\le C_0||U||,\quad\forall\; \lambda\in\mathbb R, \;|\lambda|\geq\lambda_0.
\end{align} \\
Thus, let
$\lambda\in\mathbb R$ with $|\lambda|\geq\lambda_0$, $U=(f,g,h,\ell)\in {\cal H}$, and let $Z=(u,v,w,z)\in
D({\cal A}_{\ta})$ such that
\begin{align}\label{dif2}
(i\lambda-{\cal A}_{\ta})Z=U.
\end{align}
Multiply
both sides of \eqref{dif2} by $Z$, then take the real part of the inner
product in ${\cal H}$  to
derive:
\begin{align}\label{dif3}
\gamma|A^{\theta\over2}(v+z)|^2={Re}(U,Z)\leq
||U||||Z||.
\end{align}
Equation \eqref{dif2} may be
rewritten as:
\begin{equation}\label{dif4}
\begin{cases}
i\lambda u-v=f,\\
i\lambda v+aAu+\gamma A^{\theta}(v+z)= g,\\
i\lambda w-z=h,\\
i\lambda z+bAw+\gamma A^{\theta}(v+z)=\ell.
\end{cases}
\end{equation}
Thus, \eqref{dif1} will be established as soon as
we prove the following estimate:
\begin{align}\label{dif5}
 \log(|\lambda|) ||Z||\leq C_0 ||U||,\quad\forall \lambda\in{\mathbb R}\mbox{ with } |\lambda|\geq\lambda_0.
 \end{align}
 First, we will estimate $|\lambda|||v+z||_{-1}$. Then using an interpolation inequality, get an estimate for $\log(|\lambda|)|v+z|$. Afterwards, we shall estimate $-2(\log(|\lambda|))^2\Re(v,z)$; the latter two estimates then yield estimates for both $|v|$ and $|z|$. Finally the latter estimates will be used to estimate $||u||_1$ and $||w||_1$, thereby completing the proof of \eqref{dif5}.\\
$\underline{\text {\sl Estimate }|\lambda|||v+z||_{-1}}$. 

It follows of \eqref{dif4}:
$$i\lambda( v+z)=-aAu-bAw-2\gamma A^{\theta}(v+z)+ g+\ell.$$
Therefore

$$|\lambda|||v+z||_{-1}\leq C_0(||u||_1+||w||_1+|A^{\theta-{1\over2}}(v+z)|+||U||).$$
Now, we have the following estimate
\begin{equation}\label{dif6}\begin{array}{ll}
|A^{\theta-{1\over2}}(v+z)|&\leq C_0||v+z||_{-1}^{1-\theta\over1+\theta}|A^{\theta\over2}(v+z)|^{2\theta\over1+\theta}\\&\leq C_0||(v+z)||_{-1}^{1-\theta\over1+\theta}(||U||||Z||)^{\theta\over1+\theta},\text{ thanks to }\eqref{dif3}. \end{array}\end{equation}
Consequently
\begin{equation}\label{dif7}
|\lambda|||v+z||_{-1}\leq C_0(||Z||+|\lambda|^{-{1-\theta\over1+\theta}}||\lambda(v+z)||_{-1}^{1-\theta\over1+\theta}(||U||||Z||)^{\theta\over1+\theta}+||U||).
\end{equation}
Using Young inequality, we then derive
\begin{equation}\label{dif70}
|\lambda|||v+z||_{-1}\leq C_0(||Z||+|\lambda|^{-{1-\theta\over2\theta}}(||U||||Z||)^{1\over2}+||U||).
\end{equation}
Reporting \eqref{dif70} in \eqref{dif6}, we find
\begin{equation}\label{dif60}\begin{array}{ll}
|A^{\theta-{1\over2}}(v+z)|&\leq C_0|\lambda|^{-{1-\theta\over1+\theta}}(||Z||+|\lambda|^{-{1-\theta\over2\theta}}(||U||||Z||)^{1\over2}+||U||)^{1-\theta\over1+\theta}(||U||||Z||)^{\theta\over1+\theta}. \end{array}\end{equation}
$\underline{\text {\sl Estimate }\log(|\lambda|)|v+z|}$. To this end, we rely on an interpolation inequality, (keeping in mind \eqref{log} and calling upon \eqref{dif70}):
\begin{equation}\label{dif8}\begin{array}{ll}
\log(|\lambda|)|v+z|&\leq C_0\log(|\lambda|)||v+z||_{-1}^{\theta\over1+\theta}||v+z||_\theta^{1\over1+\theta}\\&\leq \frac{C_0\log(|\lambda|)}{|\lambda|^{\theta\over1+\theta}}||\lambda(v+z)||_{-1}^{\theta\over1+\theta}||v+z||_\theta^{1\over1+\theta}\\&\leq 
 \frac{C_0\log(|\lambda|)}{|\lambda|^{\theta\over1+\theta}}
 \left(||Z||+|\lambda|^{-{1-\theta\over2\theta}}(||U||||Z||)^{1\over2}+||U||\right)^{\theta\over1+\theta}
 (||U||||Z||)^{1\over2+2\theta}\\&\leq C_0\left(||Z||^{2\theta+1\over2\theta+2}||U||^{1\over2\theta+2}+
 ||U||^{1\over2}||Z||^{1\over2}+||U||^{2\theta+1\over2\theta+2}||Z||^{1\over2\theta+2}\right). \end{array}\end{equation}
 
 $\underline{\text {\sl Estimate }-2(\log(|\lambda|))^2\Re(v,z)}$. Taking the inner product of both sides of the second equation in \eqref{dif4} and $bz$, we derive
 
 \begin{equation}\label{dif9}
b\Re(v,z)=\Re{1\over i\lambda}\left(-ab(A^{1\over2}u,A^{1\over2}z)-\gamma(A^{\theta\over2}(v+z),bA^{\theta\over2}z)+b(g,z)\right).  
\end{equation}
 Proceeding similarly, but now using the fourth equation in \eqref{dif4}, and $av$, we find
 
 \begin{equation}\label{dif10}
-a\Re(v,z)=\Re{1\over i\lambda}\left(ab(A^{1\over2}w,A^{1\over2}v)+\gamma(A^{\theta\over2}(v+z),aA^{\theta\over2}v)-a(\ell,v)\right).  
 \end{equation}
 Combining \eqref{dif9} and \eqref{dif10}, then multiplying both sides by $(\log(|\lambda|))^2$, we derive, (keeping in mind $2\beta=a-b\not=0$):
 \begin{equation}\label{dif11}\begin{array}{ll}
 -2\beta(\log(|\lambda|))^2\Re(v,z)&=(\log(|\lambda|))^2\Re{ab\over i\lambda}\left((A^{1\over2}w,A^{1\over2}v)-(A^{1\over2}u,A^{1\over2}z)\right)\\&\hskip.3cm+(\log(|\lambda|))^2\Re{1\over i\lambda}\left(\gamma(A^{\theta\over2}(v+z),aA^{\theta\over2}v-bA^{\theta\over2}z)+b(g,z)-a(\ell,v)\right).  
 \end{array}\end{equation}
 Using the first and third equations of \eqref{dif4}, it follows
 \begin{equation}\label{dif12}\begin{array}{ll} 
 \left((A^{1\over2}w,A^{1\over2}v)-(A^{1\over2}u,A^{1\over2}z)\right)=\left((A^{1\over2}w,i\lambda A^{1\over2}u-A^{1\over2}f)-(A^{1\over2}u,i\lambda A^{1\over2}w-A^{1\over2}h)\right).
  \end{array}\end{equation}
 Now, notice that
  \begin{equation}\label{dif13}\begin{array}{ll}  
\Re{ab\over i\lambda}\left(  (A^{1\over2}w,i\lambda A^{1\over2}u)-(A^{1\over2}u,i\lambda A^{1\over2}w)\right)&=
\Re{ab\over i\lambda} \left( (-i\lambda)(A^{1\over2}w, A^{1\over2}u)+i\lambda (A^{1\over2}u, A^{1\over2}w)\right)\\&=
\Re{ab}  \left(-(A^{1\over2}w, A^{1\over2}u)+ (A^{1\over2}u, A^{1\over2}w)\right)\\&=0

 \end{array}\end{equation}as the complex number inside the parentheses is purely imaginary.\\
 Consequently, \eqref{dif11} reduces to
 \begin{equation}\label{dif14}\begin{array}{ll}
 -2(\log(|\lambda|))^2\Re(v,z)&=(\log(|\lambda|))^2\Re{ab\over i\beta\lambda}\left(-(A^{1\over2}w,A^{1\over2}f)+(A^{1\over2}u,A^{1\over2}h)\right)\\&\hskip.3cm+(\log(|\lambda|))^2\Re{1\over i\beta\lambda}\left(\gamma(A^{\theta\over2}(v+z),aA^{\theta\over2}v-bA^{\theta\over2}z)+b(g,z)-a(\ell,v)\right).  
 \end{array}\end{equation}from which, we derive without any particular difficulty
 \begin{equation}\label{dif15}\begin{array}{ll}
 -2(\log(|\lambda|))^2\Re(v,z)&\leq C_0\frac{(\log(|\lambda|))^2}{|\lambda|}\left(||U||||Z||+||v+z||_\theta(||v||_\theta+||z||_\theta)\right).
 \end{array}\end{equation}
 We have the interpolation inequality, (the same holds for $||z||_\theta$):
 $$||v||_\theta\leq C_0|v|^{1-\theta}||v||_1^{\theta}\leq C_0||Z||^{1-\theta}||i\lambda u-f||_1^\theta\leq C_0(|\lambda|^\theta||Z||+||Z||^{1-\theta}||U||^{\theta}).$$
 Reporting that in \eqref{dif15}, and using \eqref{dif3}, we find, (keeping \eqref{log} in mind):
  \begin{equation}\label{dif16}\begin{array}{ll}
 -2(\log(|\lambda|))^2\Re(v,z)&\leq C_0\frac{(\log(|\lambda|))^2}{|\lambda|^{1-\theta}}\left(||U||||Z||+||Z||^{3\over2}||U||^{1\over2}+||Z||^{3-2\theta\over2}||U||^{1+2\theta\over2}\right)\\
 &\leq C_0\left(||U||||Z||+||Z||^{3\over2}||U||^{1\over2}+||Z||^{3-2\theta\over2}||U||^{1+2\theta\over2}\right).
 \end{array}\end{equation}
 Squaring \eqref{dif8} and combining obtained estimate with \eqref{dif16}, we derive
  \begin{equation}\label{dif17}\begin{array}{ll}
 (\log(|\lambda|))^2(|v|^2+|z|^2)&\leq C_0\left(||U||||Z||+||Z||^{3\over2}||U||^{1\over2}+||Z||^{3-2\theta\over2}||U||^{1+2\theta\over2}\right)
 \\&\hskip.3cm+ C_0\left(||Z||^{2\theta+1\over\theta+1}||U||^{1\over\theta+1}+||U||^{2\theta+1\over\theta+1}||Z||^{1\over\theta+1}\right)
 \end{array}\end{equation}
 It remains to estimate $\log(|\lambda|)|A^{1\over2}u|$ and $\log(|\lambda|)|A^{1\over2}w|$. Since the process is the same for both estimates, we just prove one of them.\\
 $\underline{\text {\sl Estimate }(\log(|\lambda|))^2|A^{1\over2}u|^2}$. The first step consists in estimating $\lambda^{-2}(\log(|\lambda|))^2|A^{1\over2}v|^2$, then using the first equation in \eqref{dif4} to derive the desired estimate. \\ For this purpose, we note that the second equation in \eqref{dif4} may be recast as:
 $$aAv=\lambda^2v-aAf-i\lambda\gamma A^\theta(v+z)+i\lambda g.$$
 Taking the inner product of both sides of that equation and $\lambda^{-2}(\log(|\lambda|))^2v$, we obtain
 \begin{equation}\label{dif18}\begin{array}{ll} 
 a\lambda^{-2}(\log(|\lambda|))^2||A^{1\over2}v|^2&=\log(|\lambda|))^2|v|^2-a\frac{(\log(|\lambda|))^2}{\lambda^2}\Re(A^{1\over2}f,A^{1\over2}v)\\&+\frac{(\log(|\lambda|))^2}{\lambda}\Re\left[-i\gamma( A^{\theta-{1\over2}}(v+z),A^{1\over2}v)+i (g,v)\right].\end{array}\end{equation}
 We shall now estimate the last three terms in \eqref{dif18}. Applying the Cauchy-Schwarz inequality, we obtain
 the estimate
 \begin{align*}\left|-a\frac{(\log(|\lambda|))^2}{\lambda^2}\Re(A^{1\over2}f,A^{1\over2}v)\right|&\leq
 a\frac{(\log(|\lambda|))^2}{\lambda^2}|A^{1\over2}f||A^{1\over2}v|\\&\leq {a\over4}\lambda^{-2}(\log(|\lambda|))^2||A^{1\over2}v|^2+C_0|A^{1\over2}f|^2,\text{ using }\eqref{log}\\&\leq {a\over4}\lambda^{-2}(\log(|\lambda|))^2||A^{1\over2}v|^2+C_0||U||^2\end{align*}
 Similarly, one shows, invoking \eqref{dif60}
 \begin{align*}&\left|\gamma\frac{(\log(|\lambda|))^2}{\lambda}\Re(iA^{\theta-{1\over2}}(v+z),A^{1\over2}v)\right|\\&\leq
 {a\over4}\lambda^{-2}(\log(|\lambda|))^2||A^{1\over2}v|^2+C_0(\log(|\lambda|))^2|A^{\theta-{1\over2}}(v+z)|^2\\&\leq  {a\over4}\lambda^{-2}(\log(|\lambda|))^2||A^{1\over2}v|^2\\&\hskip.3cm+C_0(\log(|\lambda|))^2|\lambda|^{-{2-2\theta\over1+\theta}}(||Z||^{2\over1+\theta}||U||^{2\theta\over1+\theta}+
 ||U||||Z||+||U||^{2\over1+\theta}||Z||^{2\theta\over1+\theta})\\&\leq  {a\over4}\lambda^{-2}(\log(|\lambda|))^2||A^{1\over2}v|^2\\&\hskip.3cm+C_0(||Z||^{2\over1+\theta}||U||^{2\theta\over1+\theta}+
 ||U||||Z||+||U||^{2\over1+\theta}||Z||^{2\theta\over1+\theta}),\text{ thanks to }\eqref{log}.\end{align*}
 and, calling upon \eqref{log} once more:
  \begin{align*}
  \left|\frac{(\log(|\lambda|))^2}{\lambda}\Re i (g,v)\right|\leq C_0|g||v|\leq C_0||U||||Z||.
  \end{align*}
  Gathering those estimates, and reporting the resulting estimate in \eqref{dif18}, we find
  \begin{align}\label{dif19} 
 a\lambda^{-2}(\log(|\lambda|))^2||A^{1\over2}v|^2\leq  C_0(||U||^2+||Z||^{2\over1+\theta}||U||^{2\theta\over1+\theta}+
 ||U||||Z||+||U||^{2\over1+\theta}||Z||^{2\theta\over1+\theta}).\end{align}
 Using the first equation in \eqref{dif4}, we easily derive
 \begin{align}\label{dif20} 
 a(\log(|\lambda|))^2||A^{1\over2}u|^2\leq  C_0(||U||^2+||Z||^{2\over1+\theta}||U||^{2\theta\over1+\theta}+
 ||U||||Z||+||U||^{2\over1+\theta}||Z||^{2\theta\over1+\theta}).\end{align}
 In a similar way, one proves
  \begin{align}\label{dif21} 
 b(\log(|\lambda|))^2||A^{1\over2}w|^2\leq  C_0(||U||^2+||Z||^{2\over1+\theta}||U||^{2\theta\over1+\theta}+
 ||U||||Z||+||U||^{2\over1+\theta}||Z||^{2\theta\over1+\theta}).\end{align}
 Combining \eqref{dif17}, \eqref{dif20} and \eqref{dif21}, one easily derives 
  \begin{equation}\label{dif22}\begin{array}{ll}
 (\log(|\lambda|))^2||Z||^2&\leq C_0\left(||U||||Z||+||Z||^{3\over2}||U||^{1\over2}+||Z||^{3-2\theta\over2}||U||^{1+2\theta\over2}\right)
 \\&\hskip.3cm+ C_0\left(||Z||^{2\theta+1\over\theta+1}||U||^{1\over\theta+1}+||U||^{2\theta+1\over\theta+1}||Z||^{1\over\theta+1}\right)\\&\hskip.3cm+C_0(||U||^2+||Z||^{2\over1+\theta}||U||^{2\theta\over1+\theta}+
 ||U||||Z||+||U||^{2\over1+\theta}||Z||^{2\theta\over1+\theta}).
 \end{array}\end{equation}Finally, applying Young inequality in \eqref{dif22}, one gets the claimed estimate
 $$\log(|\lambda|)||Z||\leq C_0||U||,$$ thereby completing the proof of Theorem \ref{neg}
\qed

\section{ Proof of Theorem \ref{reg}}
We recall the following notation introduced in the last section:
$$||u||_{\tau}=|A^{\tau\over2}u|\quad\forall\tau\in[-1,1],\quad\forall u\in V.$$
This proof is divided into two parts. In the first part, we shall prove the claimed Gevrey regularity, then in the second part, we will prove the stated optimality of the resolvent.\\
{\bf Part 1: Gevrey regularity.} 
Thanks to Lemma \ref{gevc}, it will be enough to prove that there exists a positive constant $C_0$ such that:
\begin{align}\label{eqtk}
 \sup\left\{|\lambda|^s||(i\lambda I-{\cal A}_{\ta})^{-1}||_{{\cal L}({\cal H})};~ \lambda\in \mathbb R,\text{ with }|\lambda|>1\right\}\leq C_0,
\end{align}
where the constant $s$ is the one stated in the theorem; the proof of the theorem will show how we got the stated  values for $s$. As in the last section, the constant $C_0$ may vary from line to line, and $C_0$ depends on the parameters of the system, but not on the frequency variable $\lambda$.\\
 To prove \eqref{eqtk}, we will
show that there exists $C_0>0$ such that for every $U\in{\cal
H}$, one has:
\begin{align}\label{eqtl}
|\lambda|^s||(i\lambda I-{\cal A}_{\ta})^{-1}U||\le C_0||U||,\quad\forall\; \lambda\in\mathbb R,\text{ with }|\lambda|>1,
\end{align}with the exponent $s$ as stated.\\
Thus, let
$\lambda\in\mathbb R$ with $|\lambda|>1$, $U=(f,g,h,\ell)\in {\cal H}$, and let $Z=(u,v,w,z)\in
D({\cal A}_{\ta})$ such that
\begin{align}\label{eqto}
(i\lambda-{\cal A}_{\ta})Z=U.
\end{align}
Multiply
both sides of \eqref{eqto} by $Z$, then take the real part of the inner
product in ${\cal H}$  to
derive the dissipativity estimate:
\begin{align}\label{eqtp}
\gamma|A^{\theta\over2}(v+z)|^2={Re}(U,Z)\leq
||U||||Z||.
\end{align}
Equation \eqref{eqto} may be
rewritten as:
\begin{equation}\label{eaf}
\begin{cases}
i\lambda u-v=f,\\
i\lambda v+aAu+\gamma A^{\theta}(v+z)= g{\color{red},}\\
i\lambda w-z=h,\\
i\lambda z+bAw+\gamma A^{\theta}(v+z)=\ell{\color{red}.}
\end{cases}
\end{equation}
Thus, \eqref{eqtl} will be established if
we prove the following estimate:
\begin{align}\label{eqtm}
 |\lambda|^s ||Z||\leq C_0 ||U||,\quad\forall \lambda\in{\mathbb R}\mbox{ with } |\lambda|>1.
 \end{align}
 To prove \eqref{eqtm}, first, we establish the following estimates for constituents of the velocity components defined below:
 \begin{equation}\label{est0}
 |\lambda|| v_1|+|\lambda|^{1\over2}|A^{1\over2}v_1|\leq 2 \, ||U||,\quad  |\lambda|| z_1|+|\lambda|^{1\over2}|A^{1\over2}z_1|\leq 2 \, ||U||.\end{equation}
\begin{equation}\label{est1} |\lambda|^s|v_2+z_2|\leq \varepsilon|\lambda|^s||Z||+C_\v||U||,\quad \forall\v>0\end{equation} and
\begin{equation}\label{est2}|\lambda|^{2s}|\Re(v_2,z_2)|\leq \varepsilon^2|\lambda|^{2s}||Z||^2+C_\v||U||^2,\quad \forall\v>0.\end{equation}
From those estimates, we shall derive 
\begin{equation}\label{est3} |\lambda|^s(|v|+|z|)\leq \varepsilon|\lambda|^s||Z||+C_\v||U||,\quad \forall\v>0,\end{equation}then using the latter estimate we shall establish
\begin{equation}\label{est4} |\lambda|^s(|A^{1\over2}u|+|A^{1\over2}w|)\leq \varepsilon|\lambda|^s||Z||+C_\v||U||,\quad \forall\v>0.\end{equation} Finally, gathering the last two estimates, one easily derives \eqref{eqtm}, thereby completing the proof of the claimed Gevrey regularity.\\ The attentive reader may fairly wonder why in this section we have to decompose the velocity components, instead of proceeding as in the previous section, {\tt viz}, without breaking them up. We could do that, but we would only get, for every $\theta$ in $(0,1/2]$:
\begin{align}\label{xes}|\lambda|^{2\theta\over2\theta+1}||Z||\leq C_0||U||,\end{align}which yields a weaker resolvent estimate than any of the claimed values for $s$.\\ We now turn to establishing the estimates stated above. 
We shall proceed in three steps. In the first step, \eqref{est1} will be established,  in the second step, we will prove \eqref{est2}, then combine \eqref{est0}-\eqref{est2} to derive \eqref{est3}. In the third step, we shall use \eqref{est3} to prove \eqref{est4}, and finally derive the desired estimate \eqref{eqtm}.\\

\noindent
{\bf Step 1.} In this step, we are going to show that for every $\v>0$, there exists a positive constant $C_\v$, independent of  $\lambda$ such that:
\begin{align}\label{eaj}
|\lambda|^s(|v_2+z_2|)\leq \varepsilon||Z||+C_\v||U||,\quad \forall\v>0,\quad\forall\lambda\in{\bb R}\text{ with }|\lambda|>1,
\end{align}where $v_2$ and $z_2$ are constituents of the velocity components $v$ and $z$ respectively, and are given below.\\
To prove \eqref{est3}, the main idea, adapted from the work of Liu and Renardy \cite{lre} dealing with the semigroup analyticity of thermoelastic plates, is to use an appropriate decomposition of each velocity component. Therefore, set
$$v=v_1+v_2, \text{ and } z=z_1+z_2$$ with
\begin{equation}\label{vel1} i\lambda v_1+Av_1=g,\quad i\lambda v_2=-aAu-\gamma A^{\theta}(v+z)+Av_1\end{equation}
and 
\begin{equation}\label{vel2} i\lambda z_1+Az_1=\ell,\quad i\lambda z_2=-bAw-\gamma A^{\theta}(v+z)+Az_1\end{equation}
One easily proves \eqref{est0} just by taking the appropriate inner products, then estimate real and imaginary parts separately. Next, we shall estimate $|\lambda|||v_2+z_2||_{-1}$. For this purpose, notice that \eqref{vel1} and \eqref{vel2} yield
$$ i\lambda( v_2+z_2)=-aAu-2\gamma A^{\theta}(v+z)+Av_1-bAw+Az_1.$$
It then follows, (keeping in mind that $\theta\leq1/2$):
\begin{equation}\label{vel3}\begin{array}{ll}|\lambda|||v_2+z_2||_{-1}&\leq C_0(||u||_1+||w||_1+|A^{\theta-{1\over2}}(v+z)|+||v_1||_1+||z_1||_1)\\&\leq C_0(||u||_1+||w||_1+|v|+|z|+||v_1||_1+||z_1||_1)\\&\leq C_0(||Z||+|\lambda|^{-{1\over2}}||U||),\text{ thanks to } \eqref{est0}.\end{array}\end{equation}
Now, we have the interpolation inequality
$$|v_2+z_2|\leq C_0||v_2+z_2||_{-1}^{\theta\over1+\theta}||v_2+z_2||_\theta^{1\over1+\theta}.$$
From which, one derives, making use of the dissipativity estimate \eqref{eqtp}:
\begin{equation}\label{vel4}\begin{array}{ll}|\lambda|^s|v_2+z_2|&\leq C_0|\lambda|^{s-{\theta\over1+\theta}}||\lambda(v_2+z_2)||_{-1}^{\theta\over1+\theta}||v+z -v_1-z_1||_\theta^{1\over1+\theta}\\&\leq C_0|\lambda|^{s-{\theta\over1+\theta}}(||Z||+|\lambda|^{-{1\over2}}||U||)^{\theta\over1+\theta}(||v+z||_\theta+||v_1||_\theta+||z_1||_\theta)^{1\over1+\theta}\\&\leq C_0|\lambda|^{s-{\theta\over1+\theta}}(||Z||+|\lambda|^{-{1\over2}}||U||)^{\theta\over1+\theta}(||U||^{1\over2}||Z||^{1\over2}+||v_1||_\theta+||z_1||_\theta)^{1\over1+\theta}.
\end{array}\end{equation}
We also have the interpolation inequality
$$||v_1||_\theta\leq C_0|v_1|^{1-\theta}||v_1||_1^\theta$$so that invoking \eqref{est0}, we obtain
\begin{equation}\label{vel5}||v_1||_\theta\leq C_0|\lambda|^{\theta-1-{\theta\over2}}||U||\leq C_0|\lambda|^{\theta-2\over2}||U||.\end{equation}
A similar estimate holds for $||z_1||_\theta.$\\ Therefore \eqref{vel4} becomes:
\begin{equation}\label{vel6}\begin{array}{ll}|\lambda|^s|v_2+z_2|&\leq C_0|\lambda|^{s-{\theta\over1+\theta}}(||Z||+|\lambda|^{-{1\over2}}||U||)^{\theta\over1+\theta}(||U||^{1\over2}||Z||^{1\over2}+|\lambda|^{\theta-2\over2}||U||)^{1\over1+\theta}\\&\leq C_0|\lambda|^{s-{\theta\over1+\theta}}(||Z||^{2\theta+1\over2+2\theta}||U||^{1\over2+2\theta}+|\lambda|^{-{\theta\over2+2\theta}}||U||^{2\theta+1\over2+2\theta}||Z||^{1\over2+2\theta})\\&\hskip.3cm+ C_0|\lambda|^{s-{\theta\over1+\theta}}(|\lambda|^{\theta-2\over2+2\theta}||U||^{1\over1+\theta}||Z||^{\theta\over1+\theta}+|\lambda|^{-2\over2+2\theta}||U||).\end{array}\end{equation}Now, we are going to replace $Z$ with $|\lambda|^sZ$, then find the best value possible for $s$; this is where the stated values for $s$ in the theorem are obtained. Doing that in \eqref{vel6} leads to
 \begin{equation}\label{vel7}\begin{array}{ll}|\lambda|^s|v_2+z_2|&\leq C_0|\lambda|^{s-{\theta\over1+\theta}-{(2\theta+1)s\over2+2\theta}}|||\lambda|^sZ||^{2\theta+1\over2+2\theta}||U||^{1\over2+2\theta}\\&\hskip.3cm+C_0|\lambda|^{-{(\theta+s)\over2+2\theta}+s-{\theta\over1+\theta}}||U||^{2\theta+1\over2+2\theta}|||\lambda|^sZ||^{1\over2+2\theta}\\&\hskip.3cm+ C_0\left(|\lambda|^{s-{\theta(1+s)\over1+\theta}+{\theta-2\over2+2\theta}}||U||^{1\over1+\theta}|||\lambda|^sZ||^{\theta\over1+\theta}+|\lambda|^{{-2\over2+2\theta}+s-{\theta\over1+\theta}}||U||\right).\end{array}\end{equation}
 At this stage, we draw the reader's attention to the fact that we need all exponents of $|\lambda|$ to be less than or equal to zero; in the last two terms, this is true for any $s$ in $[0,1]$. However, one readily checks that for the first term, we need $s\leq2\theta$, and for the second term, $s\leq3\theta/(2\theta+1)$. Now, we have $2\theta\leq3\theta/(2\theta+1)$, for $\theta$ in $(0,1/4]$, and for $\theta$ in $(1/4,1/2]$, the inequality is reversed. Hence the values of $s$ stated in the theorem. In the next step, we let $s$ be generic again, and find intervals for $s$, which we then compare to those found here. \\
Since  all exponents of $|\lambda|$ are  less than or equal to zero in \eqref{vel7}, it readily follows
 \begin{equation}\label{vel8}\begin{array}{ll}|\lambda|^s|v_2+z_2|&\leq C_0\left(|||\lambda|^sZ||^{2\theta+1\over2+2\theta}||U||^{1\over2+2\theta}+||U||^{2\theta+1\over2+2\theta}|||\lambda|^sZ||^{1\over2+2\theta}\right)\\&\hskip.3cm+ C_0\left(||U||^{1\over1+\theta}|||\lambda|^sZ||^{\theta\over1+\theta}+||U||\right).\end{array}\end{equation}Applying Young inequality, one then derives \eqref{est1}.\\
 {\bf Step 2.} We shall prove \eqref{est2} in this step. To this end, we recall the equations satisfied by $v_2$ and $z_2$:
 \begin{equation*} i\lambda v_2=-aAu-\gamma A^{\theta}(v+z)+Av_1,\quad  i\lambda z_2=-bAw-\gamma A^{\theta}(v+z)+Az_1\end{equation*}Taking the inner product of  the $v_2$-equation and $bz_2$, and the inner product of the $z_2$-equation and $-av_2$, then taking real parts, we find respectively:
 \begin{equation}\label{cp1}
b\Re(v_2,z_2)=\Re{1\over i\lambda}\left(-ab(A^{1\over2}u,A^{1\over2}z_2)-\gamma(A^{\theta\over2}(v+z),bA^{\theta\over2}z_2)+b(A^{1\over2}v_1,A^{1\over2}z_2)\right) 
\end{equation}
 and 
 \begin{equation}\label{cp2}
-a\Re(v_2,z_2)=\Re{1\over i\lambda}\left(ab(A^{1\over2}w,A^{1\over2}v_2)+\gamma(A^{\theta\over2}(v+z),aA^{\theta\over2}v_2)-a(A^{1\over2}z_1,A^{1\over2}v_2)\right).  
 \end{equation}
 Combining \eqref{cp1} and \eqref{cp2}, then multiplying both sides by $|\lambda|^{2s}$, we derive, (keeping in mind $2\beta=a-b\not=0$):
 \begin{equation}\label{cp3}\begin{array}{ll}
 -2\beta|\lambda|^{2s}\Re(v_2,z_2)&=|\lambda|^{2s}\Re{ab\over i\lambda}\left((A^{1\over2}w,A^{1\over2}v_2)-(A^{1\over2}u,A^{1\over2}z_2)\right)\\&\hskip.3cm+|\lambda|^{2s}\Re{1\over i\lambda}\left(\gamma(A^{\theta\over2}(v+z),aA^{\theta\over2}v_2-bA^{\theta\over2}z_2) +b(A^{1\over2}v_1,A^{1\over2}z_2)-a(A^{1\over2}z_1,A^{1\over2}v_2)\right).  
 \end{array}\end{equation}
 We note that
  \begin{equation}\label{cp4}\begin{array}{ll} 
 \left((A^{1\over2}w,A^{1\over2}v_2)-(A^{1\over2}u,A^{1\over2}z_2)\right)=\left((A^{1\over2}w,A^{1\over2}v-A^{1\over2}v_1)-(A^{1\over2}u, A^{1\over2}z-A^{1\over2}z_1)\right).
  \end{array}\end{equation}
 Using the first and third equations of \eqref{eaf}, it follows
 \begin{equation}\label{cp5}\begin{array}{ll} 
 \left((A^{1\over2}w,A^{1\over2}v)-(A^{1\over2}u,A^{1\over2}z)\right)=\left((A^{1\over2}w,i\lambda A^{1\over2}u-A^{1\over2}f)-(A^{1\over2}u,i\lambda A^{1\over2}w-A^{1\over2}h)\right).
  \end{array}\end{equation}
 Now, notice that
  \begin{equation}\label{cp6}\begin{array}{ll}  
\Re{ab\over i\lambda}\left(  (A^{1\over2}w,i\lambda A^{1\over2}u)-(A^{1\over2}u,i\lambda A^{1\over2}w)\right)&=
\Re{ab\over i\lambda} \left( (-i\lambda)(A^{1\over2}w, A^{1\over2}u)+i\lambda (A^{1\over2}u, A^{1\over2}w)\right)\\&=
\Re{ab}  \left(-(A^{1\over2}w, A^{1\over2}u)+ (A^{1\over2}u, A^{1\over2}w)\right)\\&=0

 \end{array}\end{equation}as the complex number inside the parentheses is purely imaginary.\\
 Consequently, \eqref{cp3} reduces to
 \begin{equation}\label{cp7}\begin{array}{ll}
  -2\beta|\lambda|^{2s}\Re(v_2,z_2)&=|\lambda|^{2s}\Re{ab\over i\lambda}\left((A^{1\over2}w,-A^{1\over2}f-A^{1\over2}v_1)+(A^{1\over2}u,  A^{1\over2}h+A^{1\over2}z_1)\right)\\&\hskip.3cm+|\lambda|^{2s}\Re{1\over i\lambda}\left(\gamma(A^{\theta\over2}(v+z),aA^{\theta\over2}v_2-bA^{\theta\over2}z_2)+b(A^{1\over2}v_1,A^{1\over2}z_2)-a(A^{1\over2}z_1,A^{1\over2}v_2)\right),  
 \end{array}\end{equation}from which, we derive 
 \begin{equation}\label{cp8}\begin{array}{ll}
|\lambda|^{2s} \left|2\Re(v_2,z_2)\right|&\leq C_0{|\lambda|^{2s-1}}\left(||U||||Z||+||v+z||_\theta(||v_2||_\theta+||z_2||_\theta)\right)\\&\hskip.3cm+ C_0{|\lambda|^{2s-1}}\left(|\Re \frac{1}{i}(A^{1\over2}v_1,A^{1\over2}z_2)|+|\Re \frac{1}{i}(A^{1\over2}z_1,A^{1\over2}v_2)|\right).  
 \end{array}\end{equation}
 We have the interpolation inequality, (the same holds for $||z_2||_\theta$):
 \begin{equation}\label{cp9}\begin{array}{ll}
 ||v_2||_\theta\leq C_0|v_2|^{1-\theta}||v_2||_1^{\theta}&\leq C_0(|v|+|v_1|)^{1-\theta}(||v||_1+||v_1||_1)^{\theta}\\&\leq C_0(||Z||+|\lambda|^{-1}||U||)^{1-\theta}(|\lambda|||Z||+||U||)^\theta\\&\leq
C_0(|\lambda|^{\theta}||Z||+||Z||^{1-\theta}||U||^\theta+|\lambda|^{2\theta-1}||U||^{1-\theta}||Z||^\theta+|\lambda|^{\theta-1}||U||) . \end{array}\end{equation}
Also, given that $Av_1=g-i\lambda v_1$, and $z_2=z-z_1$, we readily check, using \eqref{est0}:
 \begin{equation}\label{cp90}\begin{array}{ll}
\left|\Re{1\over i}(A^{1\over2}v_1,A^{1\over2}z_2)\right|&=\left|\Re{1\over i}(g-i\lambda v_1,z-z_1)\right|\\
&\leq (|g||z|+|g||z_1|+|\lambda||v_1||z|+|\lambda||v_1||z_1|)\\&\leq C_0(||U||||Z||+|\lambda|^{-1}||U||^2).
 \end{array}\end{equation} Similarly, one shows
  \begin{equation}\label{cp91}
\left|\Re{1\over i}(A^{1\over2}z_1,A^{1\over2}v_2)\right|\leq C_0(||U||||Z||+|\lambda|^{-1}||U||^2).
\end{equation}
 Reporting  \eqref{cp9}, \eqref{cp90} as well as \eqref{cp91} in \eqref{cp8}, and using \eqref{eqtp}, we find
  \begin{equation}\label{cp10}\begin{array}{ll}
|\lambda|^{2s} \left|\Re(v_2,z_2)\right|&\leq C_0{|\lambda|^{2s-1}}\left(||U||||Z||+|\lambda|^{-1}|||U||^2+(||U||||Z||)^{1\over2}(||v_2||_\theta+||z_2||_\theta)\right)\\&\leq C_0{|\lambda|^{2s-1}}\left(||U||||Z||+|\lambda|^{-1}|||U||^2+|\lambda|^{\theta}||Z||^{3\over2}||U||^{1\over2}+||U||^{1+2\theta\over2}||Z||^{3-2\theta\over2}\right)\\&\hskip.3cm+
 C_0{|\lambda|^{2s-1}}\left(||Z||^{1+2\theta\over2}||U||^{3-2\theta\over2}+|\lambda|^{\theta-1}||U||^{3\over2}||Z||^{1\over2}\right).
 \end{array}\end{equation}
 Proceeding as in Step 1, we replace $||Z||$ in \eqref{cp10} with $|||\lambda|^sZ||$, thereby getting
 \begin{equation}\label{cp11}\begin{array}{ll}
|\lambda|^{2s} \left|\Re(v_2,z_2)\right|&\leq C_0\left(|\lambda|^{s-1}||U|||||\lambda|^{s}Z||+|\lambda|^{2(s-1)}|||U||^2+|\lambda|^{\theta-1+{s\over2}}|||\lambda|^{s}Z||^{3\over2}||U||^{1\over2}\right)\\&\hskip.3cm+C_0\left(|\lambda|^{2s-1-{(3-2\theta)s\over2}}||U||^{1+2\theta\over2}|||\lambda|^{s}Z||^{3-2\theta\over2}+|\lambda|^{2s-1-{(1+2\theta)s\over2}}|||\lambda|^{s}Z||^{1+2\theta\over2}||U||^{3-2\theta\over2}\right)\\&\hskip.3cm+
 C_0{|\lambda|^{{3s\over2}+\theta-2}}||U||^{3\over2}|||\lambda|^{s}Z||^{1\over2}.
 \end{array}\end{equation}  
 One readily checks that, not only all exponents of $|\lambda|$ are less than or equal to zero, but also each maximal value of $s$ is larger than any of the values of $s$ obtained in Step 1. Therefore, we have 
  \begin{equation}\label{cp12}\begin{array}{ll}
|\lambda|^{2s} \left|\Re(v_2,z_2)\right|&\leq C_0\left(||U|||||\lambda|^{s}Z||+|||U||^2+|||\lambda|^{s}Z||^{3\over2}||U||^{1\over2}\right)\\&\hskip.3cm+C_0\left(||U||^{1+2\theta\over2}|||\lambda|^{s}Z||^{3-2\theta\over2}+|||\lambda|^{s}Z||^{1+2\theta\over2}||U||^{3-2\theta\over2}\right)\\&\hskip.3cm+
 C_0||U||^{3\over2}|||\lambda|^{s}Z||^{1\over2}.
 \end{array}\end{equation}  
 Squaring \eqref{vel8} and combining obtained estimate with \eqref{cp12}, we derive
  \begin{equation}\label{vel9}\begin{array}{ll}
|\lambda|^{2s}(|v_2|^2+|z_2|^2)&\leq C_0\left(|||\lambda|^sZ||^{2\theta+1\over 1+\theta}||U||^{1\over 1+\theta} +||U||^{2\theta+1\over 1+\theta}|||\lambda|^sZ||^{1\over1+\theta}\right)\\&\hskip.3cm+ C_0\left(||U||^{2\over1+\theta}|||\lambda|^sZ||^{2\theta\over1+\theta}+||U||^2\right)\\&\hskip.3cm
+C_0\left(||U|||||\lambda|^{s}Z||+|||U||^2+|||\lambda|^{s}Z||^{3\over2}||U||^{1\over2}\right)\\&\hskip.3cm+C_0\left(||U||^{1+2\theta\over2}|||\lambda|^{s}Z||^{3-2\theta\over2}+|||\lambda|^{s}Z||^{1+2\theta\over2}||U||^{3-2\theta\over2}\right)\\&\hskip.3cm+
 C_0||U||^{3\over2}|||\lambda|^{s}Z||^{1\over2}.
 \end{array}\end{equation}
 The application of Young inequality, then yields 
 
 \begin{equation}\label{vel10}
|\lambda|^{s}(|v_2|+|z_2|)\leq \v|\lambda|^{s}||Z||+C_\v||U||,\quad\forall\v>0.\end{equation}
One readily checks from \eqref{est0}:
\begin{equation}\label{vel11}
|\lambda|^{s}(|v_1|+|z_1|)\leq C_0||U||.\end{equation}
The combination of \eqref{vel10} and \eqref{vel11} then leads to the desired estimate \eqref{est3}.\\
 It remains to estimate $|\lambda|^{s}|A^{1\over2}u|$ and $|\lambda|^{s}|A^{1\over2}w|$. Since the process is the same for both estimates, we just prove one of them. That will be the object of the next and final step.\\
 
 \nt{\bf Step 3.} In this step, we shall estimate $|\lambda|^{s}|A^{1\over2}u|$. 
 First, we hall estimate $\lambda^{2s-2}|A^{1\over2}v|^2$, then use the first equation in \eqref{eaf} to derive the desired estimate. \\ For this purpose, we note that the second equation in \eqref{eaf} may be recast as:
 $$aAv=\lambda^2v-aAf-i\lambda\gamma A^\theta(v+z)+i\lambda g.$$
 Taking the inner product of both sides of that equation and $\lambda^{2s-2}v$, we obtain
 \begin{equation}\label{vel12}\begin{array}{ll} 
 a\lambda^{2s-2}|A^{1\over2}v|^2&=\lambda^{2s}|v|^2-a\lambda^{2s-2}\Re(A^{1\over2}f,A^{1\over2}v)\\&+\lambda^{2s-1}\Re\left[-i\gamma( A^{\theta-{1\over2}}(v+z),A^{1\over2}v)+i (g,v)\right].\end{array}\end{equation}
 We shall now estimate the last three terms in \eqref{vel12}. Applying the Cauchy-Schwarz inequality, we obtain
 the estimate, (keeping in mind $s\leq1$):
 \begin{align*}\left|-a\lambda^{2s-2}\Re(A^{1\over2}f,A^{1\over2}v)\right|&\leq
 a\lambda^{2s-2}|A^{1\over2}f||A^{1\over2}v|\\&\leq {a\over4}\lambda^{2s-2}|A^{1\over2}v|^2+C_0|A^{1\over2}f|^2\\&\leq {a\over4}\lambda^{2s-2}|A^{1\over2}v|^2+C_0||U||^2.\end{align*}
 Similarly, one shows, (as $\theta\leq{1\over2}$):
 \begin{align*}&\left|\gamma{\lambda}^{2s-1}\Re(iA^{\theta-{1\over2}}(v+z),A^{1\over2}v)\right|\\&\leq
 {a\over4}\lambda^{2s-2}|A^{1\over2}v|^2+C_0\lambda^{2s}|A^{\theta-{1\over2}}(v+z)|^2\\&\leq  {a\over4}\lambda^{2s-2}|A^{1\over2}v|^2+C_0\lambda^{2s}(|v|^2+|z|^2).\end{align*}
 Finally, using the Cauchy-Schwarz inequality once more, we find, (as $s\leq1$)
  \begin{align*}
  \left|{\lambda}^{2s-1}\Re i (g,v)\right|\leq C_0|\lambda|^s|g||v|\leq C_0||U||\lambda|^s|v|\leq C_0(||U||^2+\lambda^{2s}|v|^2).
  \end{align*}
  Gathering those estimates, and reporting the resulting estimate in \eqref{vel12}, we get
  \begin{equation}\label{vel13}\begin{array}{ll} 
 a\lambda^{2s-2}||A^{1\over2}v|^2&\leq  C_0(||U||^2+\lambda^{2s}(|v|^2+|z|^2)).\end{array}\end{equation}
 Using the first equation in \eqref{eaf}, we readily derive
 \begin{align}\label{disp1} 
 a\lambda^{2s}||A^{1\over2}u|^2\leq  C_0(||U||^2+\lambda^{2s}(|v|^2+|z|^2)).\end{align}
 In a similar way, one proves
  \begin{align}\label{disp2} 
 b\lambda^{2s}||A^{1\over2}w|^2\leq  C_0(||U||^2+\lambda^{2s}(|v|^2+|z|^2)).\end{align}
 Combining \eqref{disp1}, \eqref{disp2}, one finds 
  \begin{equation}\label{st1}\begin{array}{ll}
\lambda^{2s}||Z||^2&\leq C_0(||U||^2+\lambda^{2s}(|v|^2+|z|^2))
 \end{array}\end{equation}Finally, invoking \eqref{est3}, one gets the claimed estimate
 $$|\lambda|^s||Z||\leq C_0||U||,$$ thereby completing the proof of \eqref{eqtm}

\nt{\bf Part 2: Resolvent estimate optimality.} Let $\theta\in(0,1/2]$.
We are going to show that there exist a sequence of  positive real numbers  $(\lambda_n)_{n\geq1}$, and for each $n$, an element $Z_n\in {\cal D}({\cal A})$
such that for every $r\in(2\theta,1]$, one has:
\begin{align}\label{op1}
\lim_{n\to\infty}\lambda_n=\infty,\quad ||Z_n||=1, \quad \lim_{n\to\infty}\lambda_n^{-r}||(i\lambda_n-{\cal A}_\theta)Z_n||=0.
\end{align}
Indeed, if we have sequences $\lambda_n$ and $Z_n$ satisfying \eqref{op1}, then we set
\begin{align}\label{op2}
V_n=\lambda_n^{-r}(i\lambda_n-{\cal A}_\theta)Z_n,\qquad U_n=\frac{V_n}{||V_n||}.
\end{align}
 Therefore,
$||U_n||=1$ and
\begin{align}\label{op3}
\lim_{n\to\infty}\lambda_n^r||(i\lambda_n-{\cal A}_\theta)^{-1}U_n||=\lim_{n\to\infty}\frac{1}{||V_n||}=\infty,
\end{align}
which would establish the claimed result, thereby completing the proof of Theorem \ref{reg}.
 Thus, it remains to prove the existence of such sequences.\\
 For each
$n\geq1$, let $e_n$ be the eigenfunction of the operator $A$, and $\omega_n$ be its corresponding eigenvalue as in the proof of Theorem \ref{neg}. As in that proof, we seek $Z_n$ in
the form $Z_n=(a_ne_n,i\lambda_na_ne_n,c_ne_n,i\lambda_nc_ne_n)$, with $\lambda_n$ and the complex numbers
$a_n$ and $c_n$ chosen such that $Z_n$ fulfills the
desired conditions.\\ We recall that $\alpha=(a+b)/2$, $\beta=(a-b)/2$, but now,  for every $n\geq1$, we set:
\begin{align}\label{op4}
\lambda_n=\sqrt {a \omega_n}.
\end{align}
With that choice, we readily check that:
\begin{equation}\label{op5}\begin{array}{lll}
(i\lambda_n-{\cal A}_\theta)Z_n&= \begin{pmatrix}
   0\\
\left[(a\omega_n-\lambda_n^2)a_n+i\lambda_n\gamma\omega_n^\theta( a_n+ c_n)\right]e_n\\
0\\ (-\lambda_n^2+b\omega_n)c_n+i\lambda_n\gamma\omega_n^\theta( a_n+c_n))
e_n
 \end{pmatrix}\\
&= \begin{pmatrix}
   0\\
\left[ia^{1\over2}\gamma\omega_n^{\theta+{1\over2}} a_n+ia^{1\over2}\gamma\omega_n^{\theta+{1\over2}} c_n\right]e_n\\
0\\ \left[(-2\beta+ia^{1\over2}\gamma\omega_n^{\theta-{1\over2}}) c_n+ia^{1\over2}\gamma\omega_n^{\theta-{1\over2}} a_n\right]\omega_ne_n
 \end{pmatrix},\text{ by }\eqre{op4}.\\
\end{array}\end{equation}
If for each $n\geq1$, we now set
\begin{align}\label{op6}
a_n=-(1+2i\beta a^{-{1\over2}}\gamma^{-1}\omega_n^{{1\over2}-\theta})c_n.\end{align}
 It then follows from \eqref{op5}:

 \begin{equation}\label{op7}\begin{array}{lll}
(i\lambda_n-{\cal A})Z_n&= \begin{pmatrix}
   0\\
2\beta\omega_nc_ne_n\\0\\0
 
 \end{pmatrix}\end{array}\end{equation}
Therefore, we have:
\begin{align}\label{op8}
&\lim_{n\to\infty}\lambda_n^{-2r}||(i\lambda_n-{\cal A}) Z_n||^2\notag\\
&=\lim_{n\to\infty}4\beta^2|c_n|^2\lambda_n^{-2r}\omega_n^2 |e_n|^2 =4\beta^2a^{-r}\lim_{n\to\infty}\omega_n^{2\theta-r}\omega_n^{2-2\theta}|c_n|^2\\
&=0,\text{ for } r>2\theta, \text{ and }\notag
\end{align}provided that the sequence $(\omega_n^{2-2\theta}|c_n|^2)$ converges to  some nonzero real number, and $||Z_n||=1$.\\ One checks that
\begin{align}\label{op9}
||Z_n||^2&=a\omega_n|a_n|^2+\lambda_n^2|a_n|^2+b\omega_n|c_n|^2+\lambda_n^2|c_n|^2
\notag\\
&=2a\omega_n|a_n|^2+(b+a)\omega_n|c_n|^2)\notag\\
&=2a(1+4\beta^2a^{-1}\gamma^{-2}\omega_n^{1-2\theta})\omega_n|c_n|^2+(b+a)\omega_n|c_n|^2)\\
&=\left((3a+b)\omega_n^{2\theta-1}+8\beta^2\gamma^{-2}\right)\omega_n^{2-2\theta}|c_n|^2,\notag\end{align} so that we might just choose the complex number $c_n$ with
\begin{align}\label{op10}|c_n|={\omega_n^{\theta-1} \over\sqrt{\left((3a+b)\omega_n^{2\theta-1}+8\beta^2\gamma^{-2}\right)}}\end{align} to get $||Z_n||=1$ as desired.

\medskip

This  completes the proof of Theorem \ref{reg}.\qed

\begin{rem} Some further comments about our regularity results are in order. In Theorem \re{neg}, the estimate \eqref{noan} is valid for all $\theta$ in $(1/2,1]$ and $r$ in $(2(1-\theta),1]$. Notice that we do not know what happens if $r$ lies in $(0,2(1-\theta)]$. One may then fairly wonder whether some Gevrey regularity may be established for the same range of $\theta$; but this is by now an open problem.\\ In Theorem \re{reg}, we show that for $\theta$ in $(0,1/2]$, the semigroup $(S_\theta(t))_{t\geq0}$ is of certain Gevrey classes, depending on whether $\theta$ lies in $(0,1/4]$, or else 
$(1/4,1/2].$ Estimate \eqref{optgev} proves that the resolvent estimate leading to the Gevrey class corresponding to $\theta$ in $(0,1/4]$ is optimal. It is unknown to the authors of this contribution whether the resolvent estimate established in the case where $\theta$ lies in $(1/4,1/2]$ is optimal or not. However, estimate \eqref{optgev}, which is valid for all $\theta$ in $(0,1/2)$ seems to suggest that  the resolvent estimate obtained for $\theta$ lying in $(1/4,1/2]$ is not optimal. Therefore, more work is needed in the latter case. Also, in this contribution, we use the same operator in both equations; this means that the case involving different operators is still open. \end{rem}

\section{Proof of Theorem \ref{stab}}
The proof is divided into two parts: in the first part, the claimed exponential and polynomial decay estimates will be established. The second part will be devoted to proving the optimality of the polynomial decay estimate. \\

\nt
{\bf Part 1: Decay estimates.} The notations being as the proof of Theorem \ref{reg}, and given that we have \begin{align}\label{res0}i\bb R\subset\rho({\cal A}_\theta),\end{align} it remains to prove the following two resolvent estimates:
\begin{align}\label{res1}
\forall\theta\in[0,1]: \sup\left\{||(i\lambda I-{\cal A}_{\ta})^{-1}||_{{\cal L}({\cal H})};~ \lambda\in \mathbb R\right\}\leq C_0
\end{align} 
and 
\begin{align}\label{res2}
\forall\theta\in[-1,0): \sup\left\{|\lambda|^{2\theta}||(i\lambda I-{\cal A}_{\ta})^{-1}||_{{\cal L}({\cal H})};~ \lambda\in \mathbb R,\text{ with }|\lambda|>1\right\}\leq C_0,
\end{align}where in each estimate, and in the sequel, the constant $C_0$ depends on the parameters of the system, but never on the frequency variable $\lambda$.\\ We note that \eqref{res0} and \eqref{res1} yield the desired exponential decay estimate thanks to  \cite[Theorem 3]{hf} or  \cite[Corollary 4]{pr}, while \eqref{res0} and \eqref{res2} yield the claimed polynomial decay estimate thanks to  \cite[Theorem 2.4]{bort}. It remains to prove \eqref{res1} and \eqref{res2}. The proof of those two estimates are similar, but the proof of \eqref{res2} is technically more involved; so we shall prove only \eqref{res2}. However, we will be careful in our estimates, so that the proof of \eqref{res1} can easily follow along our proof.\\
 To prove \eqref{res2}, we will
show that for every $\theta$ in $[-1,0)$, there exists $C_0>0$ such that for every $U\in{\cal
H}$, one has:
\begin{align}\label{res3}
|(i\lambda I-{\cal A}_{\ta})^{-1}U||\le C_0|\lambda|^{-2\theta}||U||,\quad\forall\; \lambda\in\mathbb R,\text{ with }|\lambda|>1.
\end{align}\\
Thus, let $\theta$ in $[-1,0)$, and let
$\lambda\in\mathbb R$ with $|\lambda|>1$. Also let $U=(f,g,h,\ell)\in {\cal H}$, and $Z=(u,v,w,z)\in
D({\cal A}_{\ta})$ such that
\begin{align}\label{res4}
(i\lambda-{\cal A}_{\ta})Z=U.
\end{align}Thus, \eqref{res3} will be established if
we prove the following estimate:
\begin{align}\label{res7}
  ||Z||\leq C_0|\lambda|^{-2\theta} ||U||,\quad\forall \lambda\in{\mathbb R}\mbox{ with } |\lambda|>1.
 \end{align}
Multiply
both sides of \eqref{res4} by $Z$, then take the real part of the inner
product in ${\cal H}$  to
derive the dissipativity estimate:
\begin{align}\label{res5}
\gamma|A^{\theta\over2}(v+z)|^2={\Re}(U,Z)\leq
||U||||Z||.
\end{align}
Equation \eqref{res4} may be
rewritten as:
\begin{equation}\label{res6}
\begin{cases}
i\lambda u-v=f,\\
i\lambda v+aAu+\gamma A^{\theta}(v+z)= g,\\
i\lambda w-z=h,\\
i\lambda z+bAw+\gamma A^{\theta}(v+z)=\ell.
\end{cases}
\end{equation}
Using the first and third equations in \eqref{res6}, it follows from the dissipativity estimate
\begin{align}\label{res50}
\lambda^2|A^{\theta\over2}(u+w)|^2\leq C_0
(||U||||Z||+||U||^2).
\end{align}
To prove \eqref{res7}, we note that using first and third equations in \eqref{res6}, one derives the two equations
\begin{equation}\label{res8}
\begin{cases}
aAu=\lambda^2u-\gamma A^{\theta}(v+z)+ g+i\lambda f,\\
bAw=\lambda^2w-\gamma A^{\theta}(v+z)+\ell+i\lambda h.
\end{cases}
\end{equation}
Taking the inner product of the first equation in \eqref{res8} and $u$, and doing the same for the second equation and $w$, we find
\begin{equation}\label{res9}
a||u||_1^2=\lambda^2|u|^2-\gamma( A^{\theta-{1\over2}}(v+z),A^{1\over2}u)+( g+i\lambda f,u)
\end{equation}
and 
\begin{equation}\label{res10}
b||w||_1^2=\lambda^2|w|^2-\gamma( A^{\theta-{1\over2}}(v+z),A^{1\over2}w)+( \ell+i\lambda h,w).
\end{equation}
The application of the Cauchy-Schwarz inequality readily yields, (keeping in mind that $\theta-{1\over2}\leq{\theta\over2}$):
\begin{equation}\label{res11}
\left|-\gamma( A^{\theta-{1\over2}}(v+z),A^{1\over2}u)\right|\leq C_0|A^{\theta\over2}(v+z)||A^{1\over2}u|\leq C_0||U||^{1\over2}||Z||^{3\over2},\text{ by }\eqref{res5},
\end{equation}
\begin{equation}\label{res12}
|(g,u)|\leq |g||u| \leq C_0||U||||Z||\end{equation}
and 
\begin{equation}\label{res13}
|(i\lambda f,u)|\leq |f||\lambda||u|\leq C_0||U||^2+\lambda^2|u|^2.\end{equation}
Reporting those estimates in \eqref{res9}, we find
\begin{equation}\label{res14}
a||u||_1^2\leq 2\lambda^2|u|^2+C_0(||U||^{1\over2}||Z||^{3\over2}+||U||||Z||+||U||^2).\end{equation}
Similarly, one shows
\begin{equation}\label{res15}
b||w||_1^2\leq 2\lambda^2|w|^2+C_0(||U||^{1\over2}||Z||^{3\over2}+||U||||Z||+||U||^2).\end{equation}
At this stage, we want to draw the reader's attention to the fact that in the proof of \eqref{res11}, we have relaxed the range of $\theta$, so that our estimate is valid for every $\theta$ in $[-1,1]$.\\
To complete our proof of the resolvent estimate, it remains to estimate $\lambda^2(|u|^2+|w|^2)$. To this end, first, we shall estimate $\lambda^2|u+w|^2$, then $\lambda^2|\Re(u,w)|$, and combine those two estimates to get the desired estimate.\\
$\underline{\text{ Estimate of  }\lambda^2|u+w|^2}$. We have the interpolation inequality, (keeping in mind that $\theta<0$):

$$|u+w|\leq C_0|A^{\theta\over2}(u+w)|^{1\over1-\theta}|A^{1\over2}(u+w)|^{-{\theta\over1-\theta}},$$ from which we derive the following inequality
\begin{equation}\label{res16}
\begin{array}{ll}
\lambda^2|u+w|^2&\leq C_0\lambda^2|A^{\theta\over2}(u+w)|^{2\over1-\theta}|A^{1\over2}(u+w)|^{-{2\theta\over1-\theta}}\\&\leq C_0|\lambda|^{-2\theta\over1-\theta}(||U||^{1\over2}||Z||^{1\over2}+||U||)^{2\over1-\theta}||Z||^{-{2\theta\over1-\theta}},\text{ thanks to }\eqref{res50}\\&\leq 
C_0|\lambda|^{-2\theta\over1-\theta}(||U||^{1\over1-\theta}||Z||^{1-2\theta\over1-\theta}+||U||^{2\over1-\theta}||Z||^{-{2\theta\over1-\theta}}).
\end{array}
\end{equation}
$\underline{\text{ Estimate of  }\lambda^2\Re(u,w)}$. For this purpose, first, take the inner product of the first equation of \eqref{res8} and $bw$, next, do the same for the other equation and $-au$, and take real parts to get
$$ab\Re(A^{1\over2}u,A^{1\over2}w)=b\lambda^2\Re(u,w)-\gamma\Re(A^{\theta-{1\over2}}(v+z),bA^{1\over2}w)+\Re(g+i\lambda f,bw),$$
and
$$-ab\Re(A^{1\over2}w,A^{1\over2}u)=-a\lambda^2\Re(w,u)+\gamma\Re(A^{\theta-{1\over2}}(v+z),aA^{1\over2}u)-
\Re(\ell+i\lambda h,au).$$
Adding the two equations, we find (recalling $2\beta=a-b$):
\begin{equation}\label{res17}
2\beta\lambda^2 \Re(u,w)=\gamma\Re(A^{\theta-{1\over2}}(v+z),aA^{1\over2}u-bA^{1\over2}w)+\Re(g+i\lambda f,bw)
-\Re(\ell+i\lambda h,au).\end{equation}
Thanks to the first and third equations in \eqref{res6}, one readily checks
$$(i\lambda f,bw)=(f,-bi\lambda w)=(f,-bz-bh),\text{ and } (i\lambda h,-au)=(h,ai\lambda u)=(h,av+af).$$
Consequently
$$2\beta\lambda^2 \Re(u,w)=\gamma\Re(A^{\theta-{1\over2}}(v+z),aA^{1\over2}u-bA^{1\over2}w)+\Re\left((g,bw) -(\ell, au)-(f,b(z+h))+(h,a(v+f))\right).$$
The application of the Cauchy-Schwarz inequality then leads to
$$|\gamma\Re(A^{\theta-{1\over2}}(v+z),aA^{1\over2}u-bA^{1\over2}w)|\leq C_0|A^{\theta\over2}(v+z)|(|A^{1\over2}u|+|A^{1\over2}w|) \leq C_0||U||^{1\over2}||Z||^{3\over2},\text{ by }\eqref{res5}$$
and $$|(g,bw) -(\ell, au)-(f,b(z+h))+(h,a(v+f))|\leq C_0(||U||||Z||+||U||^2).$$
Hence
\begin{equation}\label{res18}
\lambda^2 |\Re(u,w)|\leq C_0(||U||^{1\over2}||Z||^{3\over2}+||U||||Z||+||U||^2).\end{equation}
The combination of \eqref{res16} and \eqref{res18} then yields
\begin{equation}\label{res19}
\begin{array}{ll}
\lambda^2(|u|^2+|w|^2)&\leq 
C_0|\lambda|^{-2\theta\over1-\theta}(||U||^{1\over1-\theta}||Z||^{1-2\theta\over1-\theta}+||U||^{2\over1-\theta}||Z||^{-{2\theta\over1-\theta}})\\&
\hskip.3cm +C_0(||U||^{1\over2}||Z||^{3\over2}+||U||||Z||+||U||^2).
\end{array}
\end{equation}
Using the first and third equations in \eqref{res6}, we then derive
\begin{equation}\label{res20}
\begin{array}{ll}
|v|^2+|z|^2&\leq 
C_0|\lambda|^{-2\theta\over1-\theta}(||U||^{1\over1-\theta}||Z||^{1-2\theta\over1-\theta}+||U||^{2\over1-\theta}||Z||^{-{2\theta\over1-\theta}})\\&
\hskip.3cm +C_0(||U||^{1\over2}||Z||^{3\over2}+||U||||Z||+||U||^2).
\end{array}
\end{equation}
Gathering \eqref{res14}, \eqref{res15}, \eqref{res19} and \eqref{res20}, we obtain
\begin{equation}\label{res21}
\begin{array}{ll}
||Z||^2&\leq 
C_0|\lambda|^{-2\theta\over1-\theta}(||U||^{1\over1-\theta}||Z||^{1-2\theta\over1-\theta}+||U||^{2\over1-\theta}||Z||^{-{2\theta\over1-\theta}})\\&
\hskip.3cm +C_0(||U||^{1\over2}||Z||^{3\over2}+||U||||Z||+||U||^2).
\end{array}
\end{equation}Finally, applying Young inequality, and keeping in mind that $|\lambda|>1$, we get
$$||Z||^2\leq C_0|\lambda|^{-4\theta}||U||^2,$$ hence the claimed resolvent estimate. The claimed semigroup decay estimate then follows from \cite[Theorem 2.4]{bort}.\\
\nt
{\bf Part 2: Optimality of the polynomial decay estimate.} Let $\theta\in[-1,0)$. Let the sequences $(\lambda_n)$ and $(Z_n)$ be given as in the second part of Theorem \ref{reg}. Then we already have for each $n$:

$$\lim_{n\to \infty}\lambda_n=\infty,\text{ and }||Z_n||=1.$$  
Let $r\geq0$. We shall now prove

$$\lim_{n\to\infty}\lambda_n^{r}||(i\lambda_n-{\cal A}_\theta)Z_n||=0,$$ 
provided $r<-2\theta$.\\
Proceeding as in the proof of Theorem \ref{reg}, we get
\begin{align}\label{op81}
&\lim_{n\to\infty}\lambda_n^{2r}||(i\lambda_n-{\cal A}_\theta) Z_n||^2\notag\\
&=\lim_{n\to\infty}4\beta^2|c_n|^2\lambda_n^{2r}\omega_n^2|e_n|_2^2=4\beta^2a^{r}\lim_{n\to\infty}\omega_n^{2\theta+r}\omega_n^{2-2\theta}|c_n|^2\\
&=0,\text{ for } r<-2\theta, \notag
\end{align}since  the sequence $(\omega_n^{2-2\theta}|c_n|^2)$ converges to a positive real number, thanks to \eqref{op10}.\\ Consequently, setting for each n:

\begin{align}\label{op2bis}
V_n=\lambda_n^{r}(i\lambda_n-{\cal A}_\theta)Z_n,\qquad U_n=\frac{V_n}{||V_n||},
\end{align}
 it then follows,
$||U_n||=1$ and
\begin{align}\label{op3bis}
\lim_{n\to\infty}\lambda_n^{-r}||(i\lambda_n-{\cal A}_\theta)^{-1}U_n||=\lim_{n\to\infty}\frac{1}{||V_n||}=\infty.
\end{align}
This proves the claimed estimate, and the proof of Theorem \ref{stab} is complete.\qed

\section{Examples of application}\vskip.2cm

Let $\Omega$ be a bounded domain in ${\bb R}^N$ with smooth boundary $\Gamma$. Typical examples of application includes, but are not limited to
\begin{enumerate}
\item {\bf Interacting membranes.} We consider the following system
\begin{equation*}\begin{array}{ll}
&y_{tt}-a\Delta y+\gamma(-\Delta)^\mu y_t+\gamma(-\Delta)^\mu z_t=0\text{ in }\Omega\times(0,\infty),\\
&z_{tt}-b\Delta z+\gamma(-\Delta)^\mu y_t+\gamma(-\Delta)^\mu z_t=0\text{ in }\Omega\times(0,\infty),\\
&y=0,\quad z=0\text{ on }\Gamma\times(0,\infty),
\end{array}\end{equation*}
with initial conditions 
\begin{equation*}
y(x,0)=y^0(x),\quad y_t(x,0)=y^1(x),\quad z(x,0)=z^0(x),\quad z_t(x,0)=z^1(x),
\end{equation*} 
where $a,\,b$ and $\gamma$ are positive constants, with $a\not=b$.
We take $H=L^2(\Omega)$, $A=-\Delta$ with $D(A)= H^2(\Omega) \cap H^1_0 (\Omega)$. Then $A$ is a densely defined, positive unbounded operator on the Hilbert space $H$ satisfying \eqre{xq1} by Poincar\'{e} inequality, and $D(A)$ is dense in the Hilbert space $H$. Moreover, $V:=D(A^{1 \over2})=H_0^1 (\Omega)$ and the injections $V \hookrightarrow H \hookrightarrow V^\prime=H^{-1}(\Omega)$ are dense and compact.

\medskip

Thus, the above system satisfies all results in Theorems \ref{neg},\,\ref{reg},\,and \ref{stab}.

\item{\bf Interacting  plates.} Consider the system of coupled plate equations given by 
\begin{equation*}\begin{array}{ll}
&y_{tt}+a\Delta^2 y+\gamma\Delta^{2\mu} y_t+\gamma\Delta^{2\mu} z_t=0\text{ in }\Omega\times(0,\infty),\\
&z_{tt}+b\Delta^2 z+\gamma\Delta^{2\mu} y_t+\gamma\Delta^{2\mu} z_t=0\text{ in }\Omega\times(0,\infty),\\
&y=0,\quad \frac{\partial y}{\partial\nu}=0,\quad z=0,\quad \frac{\partial z}{\partial\nu}=0\text{ on }\Gamma\times(0,\infty),
\end{array}\end{equation*}
\end{enumerate}
with initial conditions 
\begin{equation*}
y(x,0)=y^0(x),\quad y_t(x,0)=y^1(x),\quad z(x,0)=z^0(x),\quad z_t(x,0)=z^1(x),
\end{equation*} 
where $a,\,b$ and $\gamma$ are positive constants, with $a\not=b$. Let the operator $A$ defined in the Hilbert space $L^2(\Omega)$ by: $A=\Delta^2$\; with $D(A)= H^4(\Omega)\cap H_0^2(\Omega)$.

\medskip

It can be proved that that $A$ is positive, densely defined operator satisfying  \eqre{xq1}. The domain of $A^{1 \over2}$ is $V:=H^2_0(\Omega)$, and moreover, the injections $H^2_0(\Omega)\hookrightarrow L^2(\Omega)\hookrightarrow H^{-2}(\Omega)$ are dense and compact.

\medskip

Hence, the system of coupled plate equations satisfies all results in Theorems \ref{neg},\,\ref{reg},\,and \ref{stab}.

\section*{Acknowledgments.} The authors thank the handling editor as well as the referee for taking the time to read their work, and for suggesting corrections that have helped to improve its presentation.

\end{document}